\documentclass[a4paper,leqno]{article}
\usepackage{amssymb,amsmath,amsthm}
\usepackage{graphicx}
\usepackage[all]{xy}
\makeindex
\usepackage[bookmarks=true,pagebackref,,colorlinks,linkcolor=blue,pdfmark]{hyperref}
\newcounter{nequa}
\newcommand{\equa}{\stepcounter{nequa}\thenequa )\ \ }
\def\e{\varepsilon}
\def\f{\varphi{}}
\def\O{\Omega}
\def\ens#1{\left\{ #1\right\}}

\def\p{{\mathbb P}^3}
\def\qq{\quad ,\quad}
\def\G{{\cal G}}

\def\ps{{\cal P}{\cal S}}
\def\e{\varepsilon}
\def\f{\varphi}
\def\pt{\p\tilde\times\p}
\def\ft{F\tilde\times F}
\def\ov{\overline}
\def\Im{\hbox{Im}}
\def\triv{{\mathbb I}}
\newtheorem*{prop}{Proposition}
\newtheorem*{lemme}{Lemma}
\newtheorem*{remarque}{Remark}
\newcommand{\demo}{\noindent{\it Proof: }}
\renewcommand{\qed}{\par\rightline{\raise 3pt\hbox{{\it q.e.d.\qquad\qquad}}}\medskip}
\def\h{\hbox{Hom}}
\begin{document}
\title{Schubert Calculus according to Schubert}
\author{Felice Ronga}

\maketitle
\begin{abstract}
    We try to understand and justify Schubert calculus the way 
    Schubert did it. This is the english, extended version of a 
	previously posted preprint math.AG/0409281.
\end{abstract}
\tableofcontents
\listoffigures
\section{Introduction}
In his famous book \cite{Schubert:1879} ``Kalk\"ul der abz\"ahlende Geometrie'',
published in 1879,  Dr. Hermann C. H. Schubert has developed a method for 
solving problems of enumerative geometry, called Schubert Calculus
today, and has applied it to a great
number of cases. This book is self-contained~: given some aptitude to
the mathematical reasoning, a little geometric intuition and a good
knowledge of the german language, one can enjoy the many enumerative
problems that are presented and solved.  

Hilbert's 15th problems asks to give a rigourous foundation to
Schubert's method. This has
been largely accomplished using intersection theory (see
\cite{Kleiman:1976},\cite{Kleiman-Laksov:1972}, 
\cite{Fulton:1984}), and most of Schubert's calculations have been
confirmed.

Our purpose is to understand and justify the very method that Schubert
has used. We will also step through his calculations in some simple cases,
in order to illustrate Schubert's way of proceeding.

Here is roughly in what Schubert's method consists. First of all, we
distinguish basic elements in the complex projective space~: points, planes, lines.
We shall represent by symbols, say $x$, $y$, conditions (in
german~: {\it Bedingungen}\/) that some geometric objects have to
satisfy; the product $x\cdot y$ of two conditions represents the
condition that $x$ and $y$ are satisfied, the sum $x+y$
represents the condition that $x$ or $y$ is satisfied. The conditions 
on the basic elements that can be expressed using other basic elements
(for example~: the lines in space that must go through a given point)
satisfy a number of formulas that can be determined rather
easily by geometric reasoning.

\begin{figure}\label{F:ellipse}
	\centering
	\includegraphics{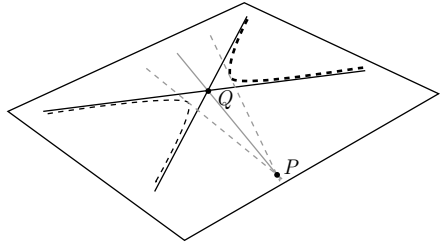}
	\caption{A degenerate solution that counts twice
	}
\end{figure}

In order to solve an enumerative problem, one tries to express it in
terms of conditions on the basic elements, by resorting if necessary
to moderatly degenerate situations, which are geometrically simpler to handle,
but might require to take in account the multiplicities of the
solutions found. For example~: to find the number of
tangents tha can be drawn from a point $P$ in a plane to a conic in
the same plane. If we take the conic to be degenerate into two distinct 
lines intersecting in a point $Q$, then the line through $P$ and $Q$ 
is the only solution, but it must be counted twice (see
figure~\ref{F:ellipse})\footnote{In the figures, the elements that
	are used to express conditions are black, the others
	are gray}.
If we degenerate the conic into a double line, then all the lines
through $P$ can be considered as tangent, and nothing can be concluded.

Schubert justifies this procedure by the Principle of the
conservation of the number ({\it Prinzip des Erhaltung der 
Anzahl}, \cite[\S\ 4, page 12]{Schubert:1879})\/, which says roughly
that that the number of solutions of an enumerative problem remains
unchanged, when its  parameters are varied, provided that this number
remains finite. In turn, Schubert justifies this principle through its algebraic
analog~: the number of solutions of a polynomial equation (in one
variable) doesn't
change if the coefficients of the polynomial are varied, provided that
multiplicities are taken into account, and that the equation doesn't
become an identity (i.e. that the polynomial is not identically
zero), in which case there are infinitely many solutions.

The strength of Schubert's approach resides in the fact that his
symbolic notation, however ambiguous, contains in germ the notion of cohomology
ring of a space (or of the Chow ring if one prefers). A condition $x$ 
represents in fact a family of conditions, that can be interpreted as 
a cohomology class in a space of configurations. It goes without
saying that when we write a product, say $x\cdot y$, the corresponding sets of
objects satisfying the conditions must be in general position (or at
least, their intersection must have the right dimension). For example,
denote by $p$ the condition, addressed to the points of $\p$, that
they must lie in a plane; $p_{g}$ denotes the condition that the points 
must lie on a line. Then the formula $p\cdot p=p_{g}$ holds. In other
words, the plane expressing the condition $p$ is generic, otherwise we
would just have $p\cdot p=p$. The ambiguity of the symbolic notation
is what makes it worthy. It must be noted that using some good
principles and some rather simple geometric reasoning, Schubert has
obtained plenty of remarkable results, whose justification, in 
accordance with today's standards of rigour, has required many great
efforts.

Among others, Schubert has
established what he has called {\it Coinzidenzformeln}\/, mainly
formula $1)$, page 44 of \cite{Schubert:1879}, which is a prototype 
of the
residual intersection formula as it can be found in 
\cite[theorem 9.2]{Fulton:1984}.
He has used this formula to establish many multiple coincidence formulas
({\it mehrfache Coinzidenzen}\/), with a rigour and effectiveness
that are not lesser than their modern analog, to be found for example in
\cite{Kazarian:2003}, although not as general.

In terms of cohomology, if $X$ denotes some space of configurations of 
geometric objects (for example, the points on a surface, the space of 
conics), a condition $x$ can be represented as the cohomology class
that is Poincar\'e dual to the fundamental homology class of a cycle $\O_{x}$ on
$X$. Then the class $x\cdot y$ is dual to $\O_{x}\cap\O_{y}$, provided
that this two cycles are in general position. The formulas proved by
Schubert on the basic elements correspond to the calculation of the
cohomology ring of the complex projective space $\p$, the grassmannian $\G$
of lines in the projective space, and finally the space $\ps$ ({\it 
Punkt un Strahl}\/), whose elements are pairs consisting of a line in 
space and a point on the line.

To give a line in $\p$ is equivalent to give a 2 dimensional vector
subspace of ${\mathbb C}^4$; with this point of view we see that there is a natural 
vector bundle $\eta$ of rank $2$ on $\G$, called the tautological bundle,
which consists of pairs $(\alpha,v)$, where $\alpha\in\G$ (regarded as 
a $2$ dimensional subspace of ${\mathbb C}^4$) and $v\in\alpha$. In
fact, the space $\ps$ in nothing but the projective bundle associated 
to $\eta$.

Note that Schubert did not introduce symbols to denote the spaces
$\p$, $\G$ and $\ps$, since somehow they constitute the ambient
universe. We shall denote by $\check\p$ the dual space of $\p$, that
is the space of projective $2$ planes in $\p$. We will assume some
familiarity with characteristic classes of vector bundles.

\section{Formulas for the basic configuration spaces}

We introduce symbols which denote geometric objects in the various
basic configuration spaces $\p$, $\check\p$, $\G$ and $\ps$. The same 
symbols will denote basic conditions imposed on the basic objects. The sets of
basic objects that are thus defined generate the homology of the
spaces; in the case of $\p$ and $\G$ they even provide a minimal cell
decomposition, that is a cell decomposition such that each
cell represents a homology class, and this homology classes are a set
of free additive generators of the homology. By expressing their 
intersections in terms of basic conditions,
the cohomology ring of these spaces will be determined explicitely.

Of course, we will use the same notation as Schubert, which is based
on the german names of the various objects. It is therefore useful to 
recall some german words~:
$$
\begin{tabular}{rl}
    Punkt :& point\\
    Gerade :& line\\
    Ebene :& plane\\
    Strahl:& litterally : ray; here it will denote most of the time
	the lines lying in a given plane\\
	&going through a given point in the plane, that is a pencil of
	lines.\\
	&Sometimes this word is synonym of line, like in {\it Punkt und
	Strahl}\\
	Fl\"ache :&surface.
\end{tabular}
$$

Note that, lacking a more precise word, we shall use {\it condition}\/
for the german  {\it Bedingung}\/ to denote a requirement imposed to
some geometric objects.

We shall work with the cohomology ring of spaces, but the Chow ring
could be used as well.

When a formula is numbered, the number is the same as in 
\cite{Schubert:1879}.

\subsection{The complex projective space $\p$}

The basic conditions that can be put on the points of space are~:
$$
\begin{tabular}{cl}
    Notation &\hfil Condition \hfil\\
    $p$ & the point must lye in a given plane\\
    $p_g$ & the point must lye on a given line\\
    $P$ & the point itself is given
    \end{tabular}
$$
The following relations are easily verified~:
$$
\equa
p^2=p_g
\qq
\equa p^3=p\cdot p_g
\qq
\equa p\cdot p_{g}=P
\qq
\equa p^3=P
\quad .
$$
As an example, the pedantic geometric interpretation of the first
formula goes as follows~: let $e_{1},\, e_{2}\subset\p$ be two
planes and
$$
\Omega_{e_i}=\{P\in\p\mid P\in e_i\}
\quad ,\quad
i=1,2
\quad .
$$
Then $p^2$ denotes the points in $\Omega_{e_1}\cap\Omega_{e_2}$ 
when $e_{1}$ and $e_{2}$ are in general position, that is when they
intersect in a line $g$; the set of points that must be on a line
has been denoted by $p_{g}$.

Let us now interpretate this formulas in cohomology. Denote by 
$t\in H^2(\p,{\mathbb Z})$ the dual class of the cycle constituted by 
the points in a plane of $\p$; then $t^2$ is the dual to the cycle
constitued by the points on a line, and $t^3$ is dual to the $0$ cycle
constituted by a single point.

If we choose a flag $p\in g\subset e$, denoting by
$\Omega_p$, $\O_g$ and $\O_e$ the corresponding sets, 
then $\O_p\subset\O_g\subset\O_e\subset\p$ is a minimal cell
decomposition of $\p$.

The case of $\check\p$, the space of $2$ planes in $\p$, 
can be treated in a similar way~:
$$
\begin{tabular}{cl}
    Notation&\hfil Condition\hfil \\
    $e$ & the plane must go through a given point\\
    $e_g$ & the plane must contain a given line\\
    $E$ & the plane itself is given
    \end{tabular}
$$
We have the formulas~:
$$
\equa e^2=e_g
\qq
\equa e^3=e\cdot e_{g}
\qq
\equa e\cdot e_g =E
\qq
\equa e^3=E
\quad .
$$
\subsection{The grassmannian ${\cal G}$ of lines in $\p$}
Here are the basic conditions~:
$$
\begin{tabular}{clc}
    Notation &\hfil Condition\hfil &Dimension \\
    $g$ & the line must cut a given line&3\\
    $g_e$ & the line must lie in a given plane&2\\
    $g_p$&the line must go through a given point&2\\
    $g_s$ & the line must belong to a given pencil&1\\
    $G$& the line itself is given&0
    \end{tabular}
$$
Choose a flag $P\in g\subset e\subset\p$ and denote by
 $\O_{g}$, $\O_{e}$, $\O_{p}$, $\O_{s}$, $\O_{G}=\ens{G}$ 
the sets of lines satisfying conditions $g$, $g_{e}$, 
$g_{p}$, $g_{s}$ and $G$ respectively. We have a diagram of inclusions~:
$$
\xymatrix{
&&\O_{p}\ar[dr]&&\\
\O_G\ar[r]&\O_{s}\ar[ur]\ar[dr]&&\O_{g}\ar[r]&{\cal G}\\
&&\O_{e}\ar[ur]&&
}
$$
and the $\O_{\bullet}$ are the cells of a minimal cell decomposition of
${\cal G}$ (see \cite[\S\ 6]{Milnor:1974}). These cells are called
{\it Schubert cycles}\/.
\begin{figure}
	\centering
	\includegraphics{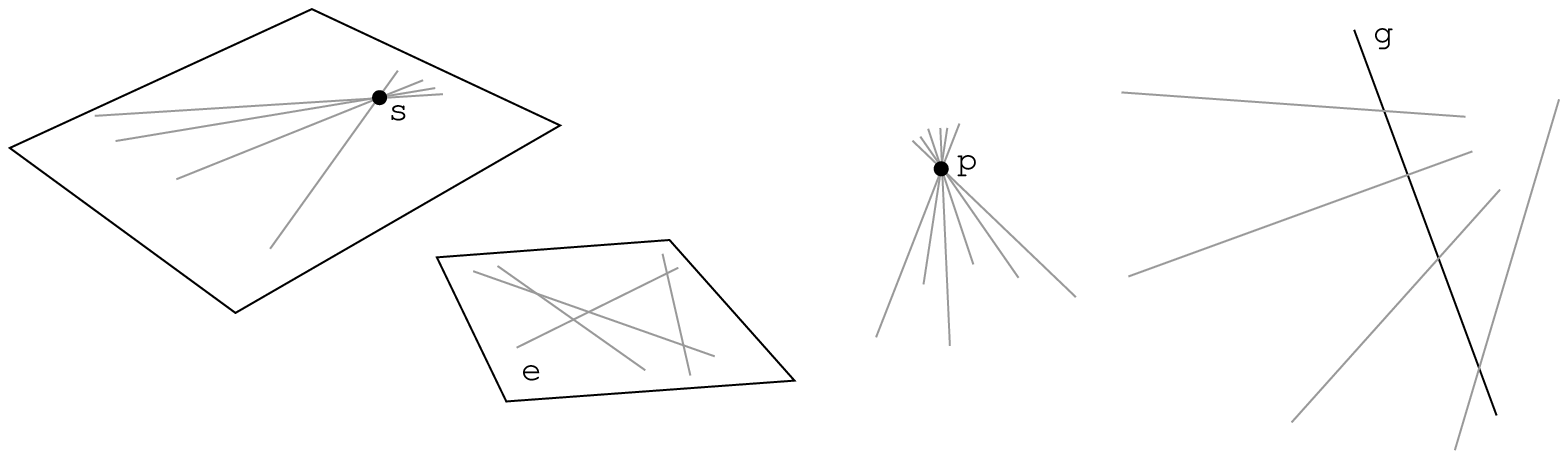}
	\caption{Schubert cycles
	\protect{$\O_{s},\,\O_{e}\,\O_{p},\,\O_{g}$}}\label{F:cycles}
\end{figure}

\begin{figure}
	\centering
	\includegraphics{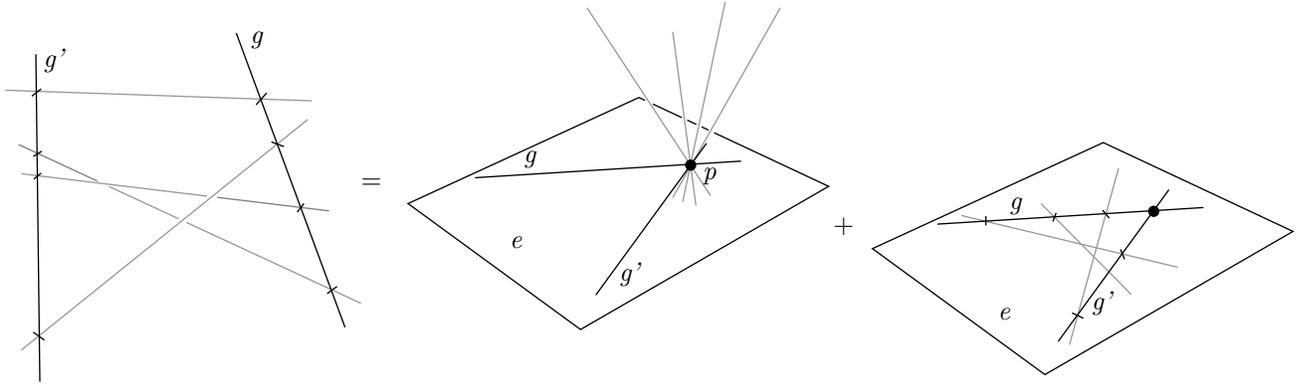}
	\caption{How to see that \protect{$g^2=g_{p}+g_{e}$}} \label{F:formula9}
\end{figure}

Now we will compute all the possible products of the basic conditions 
$g$, $g_{p}$, $g_{e}$, $g_{s}$.
In order to express $g^2$ in terms of basic conditions,
we suppose that the two given lines $g$ and $g'$  intersect in a point
$P$; by taking $e$ to be the plane of $g$ and $g'$, we have~:
$$
\O_{g}\cap\O_{g'}=\O_{p}\cup \O_{e}
$$
(see figure \ref{F:formula9}) and from this Schubert deduces, by invoking the principle of
conservation of the number, since $g$ and $g'$ are not in general
position, that~:
$$
\equa g^2=g_{p}+g_{e}\quad .
$$
We will justify this formula in two ways~: first, by expressing the
calculations in cohomology. Secondly, by showing
that $\O_{g}$ and $\O_{g'}$ intersect transversally, outside the locus 
$\O_{s}$ of lines in $e$ through $P$, which is of lower dimension;
this justifies Schubert's procedure~: in spite of the fact that the
situation is degenerated because $g$ and $g'$ are in a same plane,
the intersection of $\O_{g}$ and $\O_{g'}$  has the right dimension 2 
and there is no multiplicity to be taken into account.

Note that the linear group ${\cal G}\ell(4,{\mathbb C})$ 
acts transitively on $\G$. Therefore it can be used to put cycles in
general position~: if we choose generic lines, points or planes, the
corresponding Schubert cycles will be transversal.
It is straghtforward to check the following formulas~:
$$
\displaylines{
    \equa g\cdot g_p=g_s\qq\equa g\cdot g_e=g_s\cr
    \equa g\cdot g_s=G\qq\equa g_p\cdot g_e=0
    }
$$
By multiplying $9)$ by $g$ and using $10)$ and $11)$ we get~:
$$
\equa g^3= g\cdot g_p+g\cdot g_e
\qq
\equa g^3=2\cdot g_s
$$
By multiplying by $g$ again~:
$$
\equa g^4=2\cdot g\cdot g_s=2\cdot g^2\cdot g_e=2\cdot g^2\cdot g_p =
2\cdot g_p^2=2\cdot g_e^2=2\cdot G
\quad .
$$
Note that the formula $g^4=2\cdot G$ tells us that there are $2$ lines
that cut four given lines in general position.
This is a first and very often quoted example of application of
Schubert calculus to enumerative geometry (more on this in \S\ 2.2.2).

\subsubsection{The cohomology ring of ${\cal G}$}

Let us consider the grassmannian $\G$ as the space of 2 dimensional
vector subspaces of ${\mathbb C}^4$.
Let $\eta=(E\stackrel{\pi}{\to} \G)$ be the tautological
vector bundle of  rank $2$~:
$$
E=\ens{(\alpha,v)\in\G\times{\mathbb C}^4\mid v\in \alpha}
\qq
\pi(\alpha,v)=\alpha
\quad .
$$
Let $c_{i}(\eta)\in H^{2i}(\G,{\mathbb Z})$, $i=1,2$, be the 
Chern classes of $\eta$, and $s_{i}(\eta)\in H^{2i}(\G,{\mathbb Z})$, 
$i=1,\dots ,4$, the Segre classes (see for example
\cite{Fulton:1984}). 
They are bound by the relation~:
$$
(1+c_{1}(\eta)+c_{2}(\eta))\cdot(1+s_{1}(\eta)+s_{2}(\eta)+s_{3}(\eta)
+s_{4}(\eta))=1
\quad .
$$
Denote by $\triv^n$ the trivial bundle of rank $n$, with an
unspecified basis. Since
$\eta\subset\triv^4$, we can set $\eta'=\triv^4/\eta$, and then
$c(\eta')=s(\eta)$.

Let now $x_{1}$ et $x_{2}$ be formal variables and let
$y_{1},y_{2}\in{\mathbb Z}[x_{1},x_{2}]$ be defined by the relation~:
$$
(1+x_{1}+x_{2})\cdot(1+y_{1}+y_{2}+y_{3}+y_{4})=1
$$
which amounts to set~:
$$
y_{1}=-x_{1}
\qq
y_{2}=x_{1}^2-x_{2}
\qq
y_{3}=2x_{1}x_{2}-x_{1}^3
\qq
y_{4}=x_{1}^4-x_{2}^2+3x_{1}^2x_{2}
$$
as one can easily check. It can be shown (see \cite[proposition 
page 69]{Stong:1968}) that the ring homomorphism~:
$$
{\mathbb Z}[x_{1},x_{2}]\to H^*(\G,{\mathbb Z})
\qq
x_{i}\mapsto c_{i}(\eta)
$$
induces a ring isomorphism~:
$$
{\mathbb Z}[x_{1},x_{2}]/I(y_{3},y_{4})
\stackrel{\simeq}{\longrightarrow} H^*(\G,{\mathbb Z})
$$
 where $I(y_{3},y_{4})$ denotes the ideal generated by
 $y_{3}$ et $y_{4}$. It follows that $H^*(\G,{\mathbb Z})$ is
 generated as a group by~:
 $$
 c_{1}
 \qq 
 c_{1}^2
 \qq
 c_{2}
 \qq
 c_{1}c_{2}
 \qq
 c_{2}^2
 $$
 and the ring structure is determined by the  relations 
 $2c_{1}c_{2}-c_{1}^3=0$, $c_{1}^4-c_{2}^2+3c_{1}^2c_{2}=0$, whence
 $2c_{1}^2c_{2}-c_{1}^4=0$ and 
 $c_{1}^2c_{2}=c_{2}^2$. 
 \begin{remarque}
	 In \cite[proposition page 69]{Stong:1968}, it is asserted that 
     $H^*(\G)\simeq{\mathbb Z}[c_{1},c_{2}]/I(\ens{s_{j},j>2})$ where 
     $s_{j}$ are defined for all positive $j$ by the relations~:
     $$
     (1+c_{1}+c_{2})(1+s_{1}+s_{2}+\cdots +s_{j}+\cdots)=1
     $$
    holding in the graded ring ${\mathbb Z}[c_{1},c_{2}]$.
     But it is easy to see that $s_{j}\in I(s_{1},\dots ,s_{j-1})$,
	 and so
     $$
     I(\ens{s_{j},j>2})=I(s_{3},s_{4})\quad .
     $$
 \end{remarque}
 We will express the Poincar\'e duals of the various Schubert cells in
 terms of the Chern and Segre classes of $\eta$.
 Here are the results~:
 $$
 \begin{tabular}{|c|c|c|c|c|c|c|}\hline
     Symbolic notation & -- & $g$ & $g_{p}$ & $g_{e}$ & $g_{s}$ & $G$ \\ \hline
     Cycle & $\G$ & $\O_{g}$ & $\O_{p}$ & $\O_{e}$ & $\O_{s}$ & $\O_{G}$ \\ \hline
      Dual class & 1 &$s_{1}$ & $s_{2}$ & $c_{2}$ & $s_{1}c_{2}$ & 
      \vphantom{$\int_{0}^1$}$c_{2}^2=s_{2}^2$ \\ 
     \hline
 \end{tabular}
 $$
To do so, let      $v_{i}$, $i=1,\dots ,4$ be a basis of  ${\mathbb
C}^4$; we will denote by $\langle v_{i_{1}},\dots ,v_{i_{k}}\rangle$
the space generated by $v_{i_{1}},\dots ,v_{i_{k}}$.  The conditions
defining Schubert cycles will be expressed using the flag~:
$$
P=\langle v_{1}\rangle\subset g=\langle v_{1},v_{2}\rangle\subset
e=\langle v_{1},v_{2},v_{3}\rangle\subset{\mathbb C}^4
$$

\noindent\fbox{$\O_{g}$}

Consider the bundle morphism $\varphi_{g}:\eta\to\triv^4/\langle 
v_{1},v_{2}\rangle$ induced by the natural inclusion of $\eta$ into
$\triv^4$. Recalling that the line $g$ is the projective space
associated to the vector space $\langle v_{1},v_{2}\rangle$,
we see that 
$$
\O_{g}=\Sigma(\varphi_{g})
$$
where $\Sigma(\varphi_{g})\subset\G$ denotes the singular locus of $\varphi_{g}$,
that is the set of lines $\ell\in\G$ such that the restriction of
$\varphi_{g}$ to the fiber above $\ell$ is not injective. 
If we consider the morphism
$\Lambda^2(\varphi_{g}):\Lambda^2(\eta)\to\Lambda^2(\triv^4/\langle 
v_{1},v_{2}\rangle)$ as a section of $(\Lambda^2(\eta))^*\otimes\Lambda^2(\triv^4/\langle 
v_{1},v_{2}\rangle)\simeq (\Lambda^2(\eta))^*$, the set of zeros of
this section identifies with $\Sigma(\varphi_{g})$, and therefore its
dual class is
$c_{1} (\Lambda^2(\eta))^*=-c_{1}(\eta)=s_{1}(\eta)$.

\vspace{\medskipamount}
\noindent
\fbox{$\O_{e}$} 

Consider the natural bundle morphism $\varphi_{e}:\eta\to\triv^4/\langle 
v_{1},v_{2},v_{3}\rangle$, which corresponds to a section $\sigma$ of
$\eta^*\otimes\triv^4/\langle v_{1},v_{2},v_{3}\rangle$. Since $\O_{e}$ 
is the set of zeros of this section, its dual class is
$c_{2}(\eta^*)=c_{2}(\eta)$.

\vspace{\medskipamount}
\noindent
\fbox{$\O_{p}$} 

Here we take the natural morphism $\varphi_{p}:\langle 
v_{1}\rangle\to\triv^4/\eta$, that can be seen as a  section of 
$\triv^4/\eta$; its zeros constitute $\O_{p}$, therefore the dual
class is $c_{2}(\triv^4/\eta)=s_{2}(\eta)$.

\vspace{\medskipamount}
\noindent
\fbox{$\O_{s}$ et $\O_{G}$}

Let $e'$ be the projective plane corresponding to $\langle 
v_{1},v_{2},v_{4}\rangle$ and $g'$ the projective line corresponding
to $\langle v_{1},v_{4}\rangle$.
Notice that $\O_{s}=\O_{g'}\cap\O_{e}$ and 
$\O_{G}=\O_{e}\cap\O_{e'}$, these intersections being transversal. It 
follows that the dual classes are
$s_{1}c_{2}$ and $c_{2}^2$ respectively.

\vspace{\medskipamount}
For example, we can recover formula $9)$ by observing that
$s_{1}^2=c_{1}^2=(c_{1}^2-c_{2})+c_{2}=s_{2}+c_{2}$.
 
Also
$s_{1}^4=s_{1}(-c_{1}^3)=s_{1}(-2c_{1}c_{2})=2c_{1}^2c_{2}=2c_{2}^2$ 
shows that $g^4=2G$.

The other formulas can be recovered in a similar way.

\subsubsection{Justification of {\rm 9)} using the principle of conservation of
the number}

In order to introduce local coordinates on  $\G$, we choose a vector
subspace  $\alpha_{0}\subset{\mathbb C}^4$ of dimension $2$ and
a supplementary vector subspace $\alpha'$. Denote by
$\h(\alpha_{0},\alpha')$ the space of linear maps from
$\alpha_{0}$ to $\alpha'$. Define 
$\f:\h(\alpha_{0},\alpha')\to\G$ by associating to $A\in 
\h(\alpha_{0},\alpha')$ its graph; it is a bijection on the open
subset 
$$
U_{\alpha_{0},\alpha'}=\ens{\beta\in\G\mid \beta\cap\alpha'=\ens{0}}
\quad .
$$
It can be verified that this defines a smooth atlas on $\G$; we shall
denote by $\ell_{A}$ the projective line corresponding to $A\in 
\h(\alpha_{0},\alpha')$.

\begin{lemme}
    Let $A,B\in \h(\alpha_{0},\alpha')$ and assume that
    $\ell_{A}\in\O_{\ell_{B}}$, so that there exists a vectorial line 
    $\ell_{0}\subset\alpha_{0}$ such that
    $A|\ell_{0}=B|\ell_{0}$. 
    
    Then $\ell_{A}$ is a regular point of $\O_{\ell_{B}}$ if and only 
	if $A\not= B$, and if so~:
    $$
    T(\O_{\ell_{B}})_{\ell_{A}}
    =\ens{\ov A\in \h(\alpha_{0},\alpha')\bigm|
    \;\ov A|\ell_{0}:\ell_{0}\to\alpha'/\Im(A-B)\hbox{ is zero }}
    $$
\end{lemme}
\demo
Instead of describing $\O_{\ell_{B}}$ near $\ell_{A}$, it is easier to
work in the space
$\h(\alpha_{0},\alpha')\times \h(\ell_{0},\ell')$, where $\ell'$ is a 
vectorial line supplementary to $\ell_{0}$ in
$\alpha_{0}$. Denote by
$i_{\ell_{0}}:\ell_{0}\subset\alpha_{0}$ the inclusion, and by
$p:\h(\alpha_{0},\alpha')\times \h(\ell_{0},\ell'\to \h(\alpha_{0},\alpha'))$ 
the  projection; 
the equation
$$
(A'-B)\circ( i_{\ell_{0}}+\lambda)=0
\qq
A'\in \h(\alpha_{0},\alpha')
\, ,\,
\lambda\in \h(\lambda_{0},\lambda')
$$
defines a subset $\tilde\O$ which is in bijection through $p$ with
$\Omega_{B}\cap U_{\alpha_{0},\alpha'}$, except above 
$A'=B$. If we take the derivative of this equation at $A'=A$ we find~:
$$
\ov A\circ i_{\ell_{0}}+(A-B)\circ\ov\lambda=0
$$
where overlined symbols denote tangent vectors;
if $A\not=B$, $\hbox{Ker}(A-B)=\ell_{0}$ and so
$$
\exists\; \ov\lambda\hbox{ such that }\ov A\circ i_{\ell_{0}}+(A-B)=0
\ \Longleftrightarrow\
\ov A\circ i_{\ell_{0}}:\ell_{0}\to\alpha'/\Im(A-B)\hbox{ is zero }
$$
\qed
\begin{prop}
    Let $\ell_{B_{1}}$ and $\ell_{B_{2}}$ be two distinct lines, 
    meeting in a  point $P_{1,2}$. Then $\O_{\ell_{B_{1}}}$ and
    $\O_{\ell_{B_{2}}}$ intersect transversally except on the set
    of lines through $P_{1,2}$ that lie in the plane through $\ell_{B_{1}}$ 
    and $\ell_{B_{2}}$. 
\end{prop}
\demo
Let $\ell_{A}\in\O_{\ell_{B_{1}}}\cap\O_{\ell_{B_{2}}}$. Assume first 
that $\ell_{A}$ goes through $P_{1,2}$, and therefore is not in the
plane through $\ell_{B_{1}}$ and $\ell_{B_{2}}$. Let
$\ell_{1,2}\subset\alpha_{0}$ be the vectorial line corresponding to
$P_{1,2}$, that is such that
$B_{1}|\ell_{1,2}=B_{2}|\ell_{1,2}=A|\ell_{1,2}$. It follows from the 
lemma that
$$
T(\O_{\ell_{B_{1}}})_A\cap T(\O_{\ell_{B_{2}}})_A
=\ens{\ov A\;\bigm|\; \ov A|\ell_{1,2}:\ell_{1,2}\to 
\alpha'/\Im(A-B_{i})\hbox{ is zero },\; i=1,2}
$$
Since $A-B_{1}$ and $A-B_{2}$ have the same kernel $\ell_{1,2}$, if
they also had the same image we would have~:
$$
A-B_{1}=\lambda (A-B_{2})
$$
where $\lambda$ is a scalar, and
$\lambda\not=1$, otherwise 
$B_{1}=B_{2}$. It would follow that
$$
A=\frac{1}{1-\lambda}B_{1}-\frac{\lambda}{1-\lambda}B_{2}
$$
and $\ell_{A}$ would lie in the plane through $\ell_{B_{1}}$ and
$\ell_{B_{2}}$, a contradiction. Therefore the two  conditions 
that $\ov A|\ell_{0}\to\alpha'/\Im(A-(B_{i})$ should vanish, for
$i=1,2$,  
are independent, and therefore transversality holds.

If $\ell_{A}$ lies in the plane through $\ell_{B_{1}}$ and
$\ell_{B_{2}}$, but does not go through $P_{1,2}$, let 
$P_{1}=\ell_{A}\cap\ell_{B_{1}}$ and $P_{2}=\ell_{A}\cap\ell_{B_{2}}$
and let $\ell_{1}$, $\ell_{2}\subset\alpha_{0}$ be the vectorial lines 
such that~:
$$
(A-B_{1})|\ell_{1}=0\qq
(A-B_{2})|\ell_{2}=0
\quad .
$$
Then~:
$$
T(\O_{\ell_{B_{1}}})_A\cap T(\O_{\ell_{B_{2}}})_A
=\ens{\ov A\;\bigm|\; \ov A|\ell_{i}:\ell_{i}\to 
\alpha'/\Im(A-B_{i})\hbox{ is zero },\; i=1,2}
$$
and since $\ell_{1}\not=\ell_{2}$, these two conditions are
independent, and transversality follows.
\qed
Let $P,Q,R,S\in\p$ be $4$ points, not in a same plane, and such that $3$ 
among them are never aligned. Then, if we take the four lines
$\ell_{P,Q}$ through $P$ and $Q$, $\ell_{Q,R}$, $\ell_{R,S}$ and
$\ell_{S,P}$, one sees that the corresponding Schubert cycles
$\O_{\ell_{\bullet,\bullet}}$ intersect transversally in the two lines
$\ell_{P,R}$ and $\ell_{Q,S}$. 
Indeed, the intersection of two of the cycles is transveral according 
to the proposition, and the transversality of the remaining intersections
is elementary (for example : the intersection of the set of lines
lying in the plane through $P,Q,S$ and the set of lines in the plane
through $Q,R,S$).
\begin{figure}
	\centering
	\includegraphics{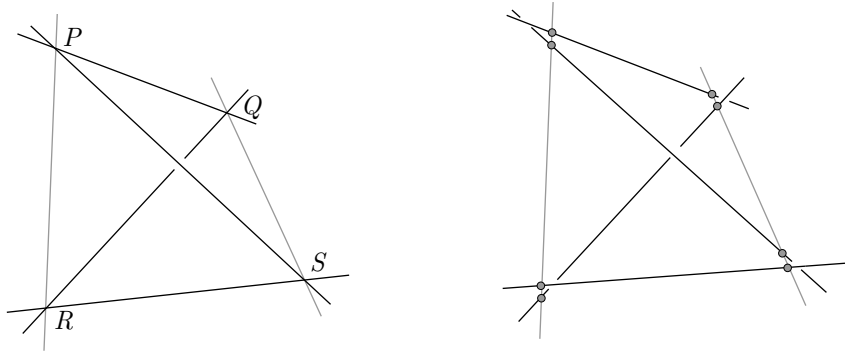}
	\caption{The $2$ lines that cut $4$ given lines}\label{F:4lines}
\end{figure}

The enumerative problem of finding the number of lines cutting $4$
given lines is often cited as an example to illustrate Schubert methods in
enumerative geometry (see \cite{Kleiman-Laksov:1972}).
What we have shown justifies to resource to the moderately degenerate case where 
each of the $4$ given lines  meets another one, and since the intersections of the
$4$ cycles are transversal, there are no multiplicities to take into account.
This last fact can be perceived with some imagination, by moving a little the 4 
black lines in figure~\ref{F:4lines}, and seeing that near each gray line
there is only one solution.

On the other hand, here is a degenerate situation that has been
pointed out to me by my colleague Alexandre Gabard. On a smooth quadric
surface in $\p$ there are two sytems of lines, which correspond to the
horizontal and the vertical lines respectively if one identifies the
quadric with ${\mathbb P}^1\times{\mathbb P}^1$.
If one takes 4 lines of one system, they might seem to be in general
position in the space ; however, any line
of the other system cuts the four given lines~: we are in a very
degenerate situation, with an infinite number of solutions.

\subsection{The space $\ps$ of  points on a line of $\p$}
Recall that the space $\ps$ is the set  of pairs consisting of a point
on a line of $\p$. In order to express conditions on elements of this
space, we shall use symbols of the form $xy$, where $x$ is y symbol
expressing a condition on points and $y$ is a symbol
expressing a condition on lines. Thus the symbol $pg$ denotes the pairs
consisting of a point on a line, where the point must lie on a given
plane, and the line must cut a given line; if we denote by
$\O_{pg}$ the set of these pairs, and yet by $g$ a line and by $e$ 
a plane, we have~:
$$
\O_{pg}=\ens{(\ell,Q)\in\ps\mid \ell\cap g\not=\emptyset\qq Q\in e}
\quad .
$$
\begin{figure}[ht]
	\centering
	\includegraphics{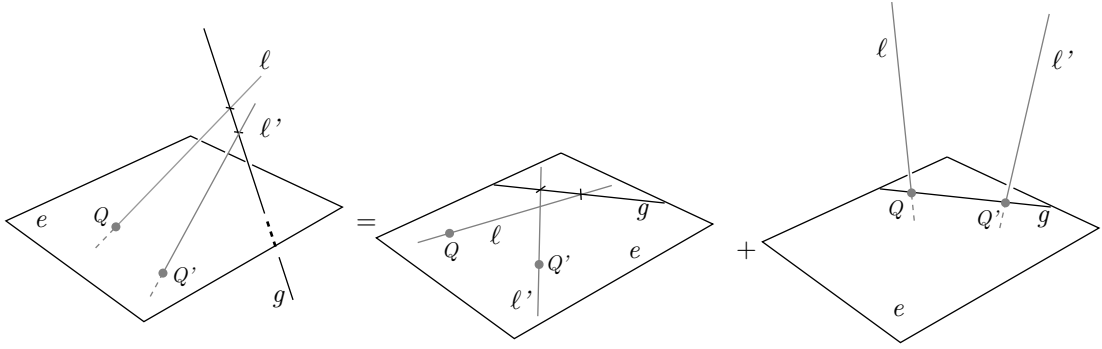}
	\caption{How to see that $pg=g_{e}+p_{g}$}\label{F:formulaI}
\end{figure}

By the principle of conservation of the number, we can take the line
$g$ in the plane $e$, in which case~:
$$
\O_{pg}=\ens{(\ell,Q)\in\ps\mid Q\in 
g}\cup\ens{(\ell,Q)\in\ps\mid\ell\subset e}
$$
(see figure \ref{F:formulaI}) and so the formula of \cite[page 25]{Schubert:1879} follows
\footnote{There is a misprint in \cite[page 25]{Schubert:1879}~: the
formula shown there is $pg=p_{g}+g$}
$$
\hbox{I})\ \ pg=p_{g}+g_{e}=p^2+g_{e}
\quad .
$$
This is a fundamental formula, in the sense that any other formula in 
$\ps$ will follow from this one and from the formulas that we have
already shown to hold in $\p$ and $\G$; the reason is explained in the
next \S\ . 

Let us show some other formulas in $\ps$ anyway.
By multiplying $\hbox{I})$ by $p$, then by 
$g$ we get~:
$$
\displaylines{
p^2g=pp_{g}+pg_{e}=p^3+pg_{e}\cr
pg_{e}+pg_{p}=pg^2=p_{g}g+g_{e}g=p_{g}g+g_{s}=p^2g+g_{s}
}
$$
and by adding the far left and far right expressions of these two
lines~:
$$
\hbox{II})\ \ pg_{p}=p^3+g_{s}
$$
and similarly one obtains (see \cite[page 26]{Schubert:1879})~:
$$
\hbox{III} )\ \ pg_{s}=p^2g_{p}=G+p^3g=G+p^2g_{e}
\quad .
$$

\paragraph{Justification of {\rm I)} using cohomology}
If we regard $\G$ as the space of 2 dimensional vector subspaces of
${\mathbb C}^4$, $\ps$ is the projective bundle associated to the
tautological bundle $\eta$ of rank $2$. The tautological line bundle
$\gamma=(F\stackrel{\pi}{\to}\ps)$ on $\ps$ can be defined as
$$
F=\ens{(\alpha,\ell,v)\in\G\times\p\times{\mathbb C}^4\mid 
\ell\subset\alpha\; ,\, v\in\ell}
\qq
\pi(\alpha,\ell,v)=(\alpha,\ell)
\quad .
$$
Set $s=c_{1}(\gamma)$. Notice that $H^*(\ps)$ is a
$H^*(\G)$-module via the homomorphism induced by the natural projection
$p:\ps\to\G$. We know that (see \cite[theorem page 62]{Stong:1968})
the ring homomorphism
$$
H^*(\G)[s]\to H^*(\ps)
\qq
s\mapsto c_{1}(\gamma)
$$
induces an isomorphism
$$
H^*(\G)[s]/I(s^2-sc_{1}(\eta)+c_{2}(\eta))\stackrel{\simeq}{\longrightarrow}
H^*(\ps)
$$
where $I(s^2-sc_{1}(\eta)+c_{2}(\eta))$ denotes the ideal generated by
the polynomial $s^2-sc_{1}(\eta)+c_{2}(\eta)$, which is nothing else
than $c_{2}(p^*(\eta)/\gamma)$ once we substitute $s$ by $c_{1}(\gamma)$;
it vanishes because $p^*(\gamma)/\gamma$ is of rank $1$. 

Let us allow to denote by the same symbol a condition on the basic
elements as well as the Poincar\'e dual to the cycle defined by this
condition. For example, using \S\ 2.2.1, we shall write
$g=s_{1}(\eta)=-c_{1}(\eta)$; we will also write $s$ for
$c_{1}(\gamma)$, so that $p=-s$.

It follows from the very definition of the symbols that
$pg=(-s)(-c_{1}(\eta))$, and by \S\ 2.2.1 $g_{e}=c_{2}(\eta)$.
On $H^*(\ps)$ we have the relation
$$
sc_{1}(\eta)=s^2+c_{2}(\eta)
$$
that can also be written
$$
pg=p^2+g_{e}
$$
which is formula I). Therefore this formula is exactly the relation by
which $H^*(\G)[s]$ has to be divided in order to obtain $H^*(\ps)$.

\section{Coincidence formulas}

If $X\subset{\mathbb P}^1\times{\mathbb P}^1$ is a curve of bi-degree
$(p,q)$, the restriction of its equation to the diagonal
is of degree $p+q$, therefore $X$ cuts this diagonal in $p+q$ 
points, counted with multiplicity. 
We can reformulate this remark by saying that $X$ is a one parameter
family of pairs $(P,Q)$ of points on the line; if there are $q$ pairs in the
family with a given first point $P$, and $p$ pairs with a given second
point $Q$, then there are $p+q$ pairs of the form $(P,P)$. This
is the {\it Principle of correspondence}\/ that has been stated and proved by
Chasles \cite[Lemme I, page 1175]{Chasles:1864}. 

We will generalize this formula according to
\cite[ pages 42 and following]{Schubert:1879}, 
using the same notation and its fruitful ambiguity.
Consider pairs of points in the projective space and imagine that when 
two points of a pair come to coincide, the line joining the two points 
has a well defined limit. Let us call $p$ and $q$ the two points of a 
pair, and $g$ the line joining them; denote by $\e$ the condition 
that $p$ and $q$ are infinitely near, but still determine the line 
joining them.

Assume that a one parameter system $X$ of such pairs of points is 
given. Note that if we still denote  by $p$ the number of pairs 
$(P,Q)\in X$ such that $P$ lies in a given plane (a condition that we 
also denoted by $p$), and by $q$ the number of pairs $(P,Q)\in 
X$ such that $Q$ lies in a given plane (a condition that we also 
denote by $q$), then $X$ is of bi-degree $(p,q)$.

\begin{figure}
	\centering
	\includegraphics{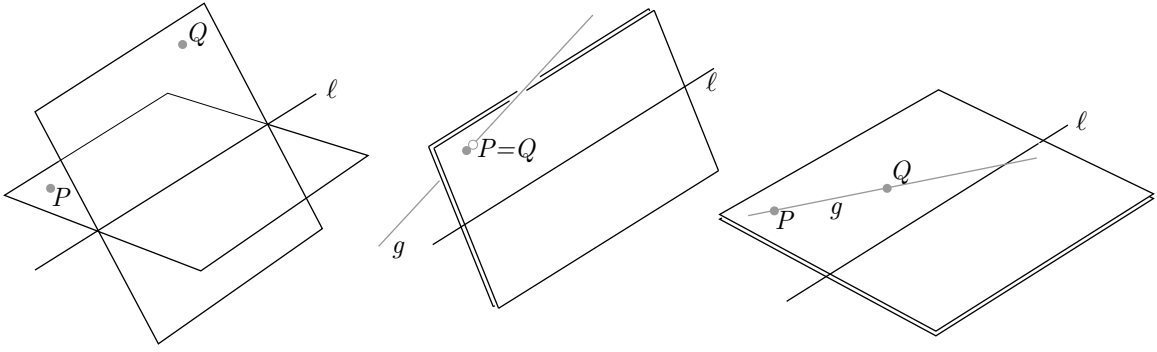}
	\caption{The space $Y$ and conditions $\e$ and 
	 $g$}\label{F:coincidence}
	\end{figure}

Now we take a line $\ell$ and consider pairs of planes through 
$\ell$, such that the first plane contains $P$, the second contains 
$Q$. These pairs constitute a curve $Y$ in the space of pairs of planes 
through $\ell$, that can be identified to  
${\mathbb P}^1\times{\mathbb P}^1$ ; $Y$ is also of
bi-degree $(p,q)$. It follows from Chasles {\it Principle of 
correspondence}\/ that there are 
$$
p+q
$$
single planes that contain a pair of points in $X$. Among these, first
we have $\e$ of them that arise because they contain a pair 
of coinciding points of $X$, secondly we have those planes who contain 
a line $g$ joining two distinct points of a pair in $X$; the latter
is equivalent to say that the line $g$ must cut the line 
$\ell$ (this kind of conditions on lines has been denoted by $g$, 
same notation as the line $g$ joining $P$ and $Q$ !); see 
figure~\ref{F:coincidence}. 
Therefore we have~
\begin{equation*}
\e=p+q-g\quad .\tag{$\spadesuit$}
\end{equation*}
This formula proves to be useful to establish enumerative 
formulas concerning special positions of lines with respect to a 
surface, among others. For example, Schubert uses it to show that
a generic surface of degree $n$ possesses
$$
\frac{1}{12}n(n-4)(n-5)(n-6)(n-7)(n^3+6n^2+7n-30)
$$
lines that are tangent at four distint points (see \cite[page 237, 
formula 21)]{Schubert:1879}.
We will follow Schubert's procedure to establish this formula in our 
third and last example of the next paragraph.

First, let us justify formula $(\spadesuit)$.
Consider the space $\p\tilde\times\p$ obtained by blowing up the 
diagonal $\Delta$ in $\p\times\p$. The map that associates to
$(P,Q)\in\p\times\p\setminus\Delta$ the line through $P$ 
and $Q$ extends to a map $\f:\pt\to\G$, such that
$$
\f^*(\Lambda^2(\eta))=\gamma^*\otimes({\cal O}(1)_{1}\otimes{\cal 
O}(1)_{2})
$$
where $\gamma$ denotes the line bundle associated to the blown up 
diagonal, ${\cal O}(1)_{i}$ the pull-back by the projection of $\pt$ 
on the i-th factor of 
$\p\times\p$ of the vector bundle  of homogeneous $1$-forms on $\p$.

To see that $\f$ has these properties, we can use 
the Pl\"ucker imbedding $\psi$ of the grassmannian $\G$ into ${\mathbb P}^5$, 
defined as follows.  If $g\in\G$, choose distinct $P,Q\in g$; if 
$P=[x_{1},\dots ,x_{4}]$, $Q=[y_{1},\dots ,y_{4}]$, et  $x=(x_{1},\dots 
,x_{4})$, $y=(y_{1},\dots ,y_{4})$, set 
$$
\psi(g)=[x\wedge y]\in{\mathbb P}(\Lambda^2({\mathbb C}^4))\simeq{\mathbb 
P}^5
\quad .
$$
It can be checked that  $\psi$ is well defined and that it is an
imbedding, whose image is
$$
\ens{[P\wedge Q]\in {\mathbb 
P}(\Lambda^2({\mathbb C}^4))\mid P\wedge Q\not=0}
\quad ,
$$ 
and this image can be identified to $\G$. Note that the pull-back
by $\psi$ of the vector
bundle ${\cal O}(1)_{{\mathbb P}(\Lambda^2({\mathbb C}^4))}$ of
homogeneous 1-forms on
${\mathbb P}(\Lambda^2({\mathbb C}^4))$ is naturally isomorphic to
$\Lambda^2(\eta^*)$.
Consider the map
$$
\Phi:{\mathbb C}^4\times{\mathbb C}^4\to\Lambda^2({\mathbb C}^4)
\qq
(x,y)\mapsto x\wedge y
\quad .
$$
Its derivative with respect to $x$, at a point point $(y,y)$, $y\not=0$, 
can be writtent $v\mapsto v\wedge y$, and its kernel is the line
supporting $y$. 
It follows from this that $\Phi$ induces a morphism
$$
\f:\pt\to\G
\hbox{\quad with\quad}
\f^*(\Lambda^2(\eta^*))\simeq
\gamma^*\otimes({\cal O}(1)_{1}\otimes{\cal 
O}(1)_{2})
\quad .
$$
With an additional effort, one can even show that
$
\f^*(\eta)=(\gamma^*\otimes{\cal O}_{1}(-1))\oplus{\cal
O}_{2}(-1)
$
but we won't use it.
Recall now from \S\ 2.2.1 that the dual class to $\Omega_{g}$ is
$s_{1}(\eta)=-c_{1}(\eta)$~; since 
$$
\f^*(-c_{1}(\eta))=
c_{1}(\gamma^*\otimes{\cal O}(1)_{1}\otimes{\cal O}(1)_{2})
$$
setting $t_{i}=c_{1}({\cal O}(1)_{i}$, $i=1,2$, 
$\e=c_{1}(\gamma)$, we get~:
$$
\f^*(-c_{1}(\eta))=t_{1}+t_{2}-\e
\quad .
$$
In order to recover Schubert's coincidence formula $(\spadesuit)$, we 
must observe that in this context the condition $g$, that is the
condition that the line through a pair $(P,Q)$ cuts a given line, corresponds to
$\f^*(s_{1}(\eta))=\f^*(-c_{1}(\eta))$; 
and the conditions $p$ and $q$ correspond to $t_{1}$ and $t_{2}$
respectively, that is the dual class to a hyperplane in the first
factor $\p$, respectively the second.

\subsection{Coincidences of intersections of a line and a surface}

We will present three examples of computations using the coincidence 
formula. The first one will be justified also using cohomology, the 
other two will be treated only the Schubert's way.

Let $F\subset\p$ be a smooth surface of degree $n$;
following \cite[page 229]{Schubert:1879}, we denote by 
$p_{1}$, $p_{2}$, \dots $p_{n}$ the points 
of intersection of a line $g$ with $F$. Denote by $\e_{2}$ the
condition that $2$ of these points coincide. Then it follows from the
coincidence formula $(\spadesuit)$ that
\begin{equation*}
\e_{2}=p_{1}+p_{2}-g\quad .\tag*{$\clubsuit$}
\end{equation*}

As usual, the same symbol $p$ is used to express a condition (to be 
in a plane) and its recipient (a point).

\subsection{First example: the class of a curve}

Let's multiply formula $\clubsuit$ by  $g_{s}$~:
$$
\e_{2}g_{s}=p_{1}g_{s}+p_{2}g_{s}-G
$$
using formula III)~:
$$
\e_{2}g_{s}=G+p_{1}^3g+G+p_{2}^3g-G=G
$$
because $p^3=0$ (the generic intersection of $3$ planes and a surface 
is empty). 
It remains to interpretate the symbol $G$ in this context~: it
represents the pairs of distinct points on the intersection of a fixed
line and the surface; there are $n(n-1)$ such pairs. Thus we recover
the formula for the class (i.e. degree of the dual) of a plane curve
of degree $n$; indeed, $\e_{2}g_{s}$ represents the lines tangent to
the surface belonging to a given pencil, which is the same as the
lines in a plane, passing through a given point, tangent to the plane curve
obtained as intersection of the surface with the plane.

\paragraph{Justification using cohomology}

Denote by $\ft$ the space obtained by blowing up the diagonal in
$F\times F$. We have  $\ft\subset\pt$, and would like to express its dual 
class. The following result will help us.

Let $X$ be a smooth variety  and $A$, $Y\subset X$ smooth
subvarieties that intersect nicely, that is such that $A\cap Y$ is
smooth, and that for all $x\in A\cap Y$~:
$$
TA_{x}\cap TY_{x}=T(A\cap Y)_{x}
\quad .
$$
Then we have an exact sequence of vector bundles~: 
$$
    0\to T(A\cap Y)\to TA|_{A\cap Y}\oplus TY|_{A\cap Y}\to
    TX|_{A\cap Y}\to E\to 0
$$
where $E$ is defined by the sequence itself;
it is called the {\it excess bundle} and $k$ 
will denote its rank. Note that $k=0$ if and only if $A$ and $Y$
intersect transversally.
\begin{prop}
    Let $X$ be a smooth variety, $A$, $Y\subset X$ smoth subvarieties 
	that intersect nicely. Then, denoting by~:  
    \vspace{\smallskipamount}
    $$
    \begin{tabular}{c c l}
    $\delta_{U,V}$&\quad& the dual class to $U$ in $V$\\
    $\tilde X$&\quad&the blowing up of $X$ along $Y$\\
    $\tilde A$&\quad& the strict transform of $A$\\
    $\e$&\quad&the dual class to the exceptional divisor in $\tilde X$\\
    $p:\tilde X\to X$&\quad& the projection of the blowing up\\
    $j:\tilde Y\subset \tilde X$&\quad&the natural inclusion
    \end{tabular}
    $$
  \vspace{\bigskipamount}\noindent  
  we have~:
    $$
    \delta_{\tilde A,\tilde X}
    =p^*(\delta_{A,X})- j_{!}\Bigg((p|\tilde Y)^*(\delta_{A\cap Y,Y})\cdot 
  \underbrace{ 
  \sum_{i=0}^{k-1}(-1)^i\e^ic_{k-i-1}(E)}_{=c_{k-1}(E/\gamma)}\Bigg)
    $$
\end{prop}
It is a special case of
\cite[Theorem 6.7]{Fulton:1984}.

As an application, consider the subvarieties
$F\times F$ and $\Delta$ of $\p\times\p$. In this case, the excess
bundle identifies to the normal bundle of $F$ in $\p$, that is
${\cal O}(n)_{\Delta}$, and $\delta_{F,\p}=nt$, therefore 
$\delta_{F\times F,\p\times \p}=nt_{1}\cdot nt_{2}=n^2t_{1}t_{2}$. It 
follows that
\begin{equation*}
\delta_{\ft,\pt}=n^2t_{1}t_{2}-nt\e\tag{$\heartsuit$}
\end{equation*}
where $t$ denotes indistinctly $t_{1}$ or $t_{2}$, since
$\e t_{1}=\e t_{2}$.

The formula $\heartsuit$ can also be proved in Schubert's spirit as follows.
Denote by
$p:\p\tilde\times\p\to\p\times\p$ the projection 
of the blowing up and set $\tilde\Delta_{F}=\tilde\Delta\cap^{-1}(F\times 
F\cap\Delta)$. Then~:
$$
p^{-1}(F\times F)=(F\tilde\times F)\cup\tilde\Delta_{F}
\quad .
$$
Taking the dual classes, we see that
$$
p^*(\delta_{F\times F,\p\times\p})
=\delta_{\ft,\pt}+\delta_{\tilde\Delta_{F},\pt}
$$
and
$$
\delta_{\tilde\Delta_{F},\pt}=
\delta_{\tilde\Delta_{F},\tilde\Delta}\cdot\delta_{\tilde\Delta,\pt}
=(nt)\cdot \e
$$
whence the formula $\heartsuit$.

In particular, taking $n=1$, i.e. $F$ is a plane, we obtain that
$$
\delta_{\f^{-1}(\O_{e})}=t_{1}t_{2}-t\e
\quad .
$$
Now $g_{s}=gg_{e}$, hence 
$\f^*(g_{s})=\f^*(g)\f^*(g_{e})=(t_{1}+t_{2}-\e)(t_{1}t_{2}-t\e)$.
In order to calculate
$\e_{2}g_{s}$ we must multiply $\e\cdot\f^*(g_{s})$ by
$\delta_{\ft,\pt}$ and evaluate this class on $\pt$, which amounts to 
evaluate $\delta_{\ft,\pt}\cdot\f^*(g_{s})$ on $\tilde \Delta$; but
$$
\langle \delta_{\ft,\pt}\cdot\f^*(g_{s}),\tilde\Delta\rangle=
\langle (n^2t_{1}t_{2}-nt\e)(t_{1}t_{2}-t\e)(t_{1}+t_{2}-\e),\tilde\Delta\rangle=
\langle (n^2t^2-nt\e)(t^2-t\e)(2t-\e),\tilde\Delta\rangle
$$
and since $t^4=0$,
$(n^2t^2-nt\e)(t^2-t\e)(2t-\e)=t^2(-n\e^3+\e^2(n^2t+3nt))$. 
Instead of evaluating on $\tilde\Delta$, we can apply the integration 
over the fibers of
$\pi$ (or Gysin homomorphism) $\pi_{!}$ and
evaluate on $\p$, where $\pi:\tilde\Delta\to\p$ is the natural
projection, that is the projection of the projective bundle associated
to $T\p$; we have the following formulas~:
$$
\pi_{!}(\e^2)=1 
\qq
\pi_{!}(\e^3)=c_{1}(T\p)=4t
\quad ,
$$
either by the very definition of Segre classe
given in \cite[\S\ 3.1]{Fulton:1984}), or using 
\cite[theorem page 62]{Stong:1968},
and so
$$
\langle t^2(-n\e+\e^2(n^2t+3nt)),\tilde\Delta\rangle
=\langle t^3(-4n+n^2+3n),\p\rangle=n(n-1)
\quad .
$$

\subsection{Second example: the number of bitangent lines to a plane curve}
This example, as well as the next and last one, will be treated 
the Schubert's way, without 
any further justification (see \cite[page 229]{Schubert:1879}).

Let $F\subset\p$ be a smooth surface and
denote by $\e_{22}$ the condition that a line is tangent at two
distinct points of $F$. This condition says that, 
among the points $p_1,\dots ,p_n$, 
intersection of the line with $F$, two pairs coincide, 
say $p_1,p_2$ and $p_3,p_4$. It follows from the coincidence 
formula $(\spadesuit)$
that~:
$$
2\cdot\e_{22}=(p_1+p_2-g)(p_3+p_4-g)
$$
where the coefficient $2$ is due to the fact that the roles of $(p_{1},p_{2})$ et $(p_{3},p_{4})$
can be exchanged on a bitangent line; and so
$$
2\cdot\e_{22}=p_1p_3+p_1p_4+p_2p_3+p_2p_4-p_1g-p_2g-p_3g-p_4g+
\underbrace{g^2}_{=g_e+g_p}
\quad .
$$
The symbols $p_ip_j$, $i\not=j$, all have the same meaning, 
and also the symbols $gp_i$; therefore we can write~:
$$
2\cdot\e_{22}=4p_1p_3-4p_1g+g_e+g_p
\quad .
$$
Now we multiply by $g_e$ : $\e_{22}g_e$ denotes the lines
bitangent to the surface, lying in a given plane; 
they are therefore the bitangent lines to the curve obtained by
intersecting the surface with the plane.
We have~:
$$
2\cdot \e_{22}g_e=4p_1p_3g_e-4p_1gg_e+g_{e}^2+ \underbrace{g_pg_e}_{=0}=
4p_1p_3g_e-4p_1g_s+G\stackrel{\hbox{(by III))}}{=}
4p_1p_3g_e-4p_1^3g-3G
\quad .
$$
\begin{figure}[ht]
	\centering
	\includegraphics{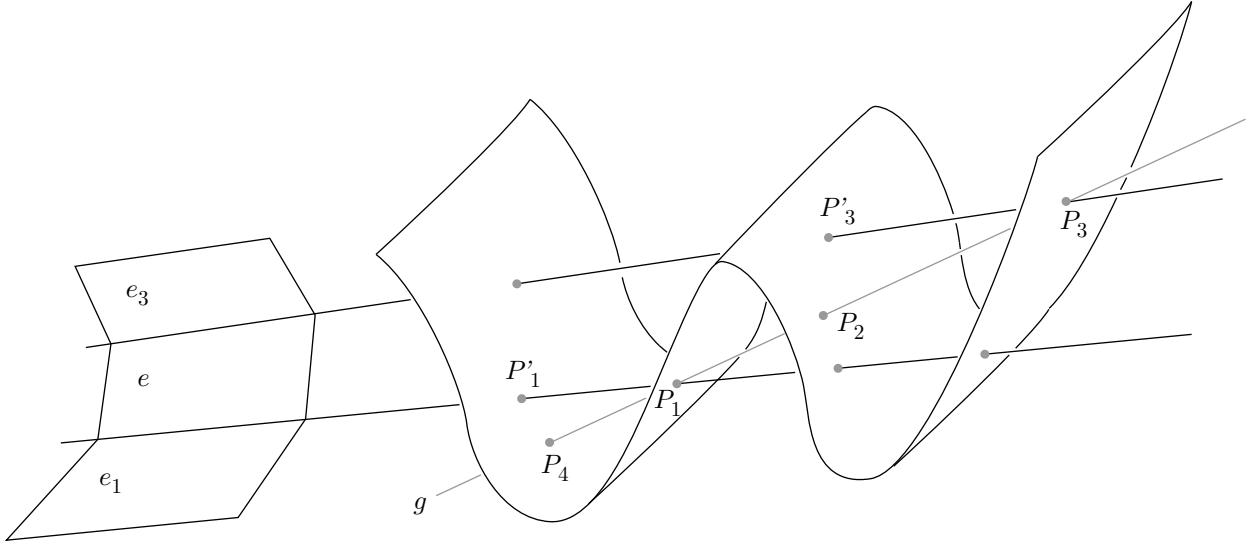}
	\caption{How to calculate $p_1p_3g_e$}
	\label{F:bitangents}
\end{figure}
Let us calculate $p_1p_3g_e$. In fact, we are working in
$$
\big(F\tilde\times F\big)\times\big(F\tilde\times F\big)\times\G
$$
and a generic element of this set can be represented by
$((P_1,P_2),(Q_1,Q_2),g)$, with $P_i,Q_i\in F\cap g$. The condition $p_i$ 
requires that $P_i$ lies in a plane $e_i$, $i=1,3$, and $g_e$ requires
that $g$ lies in a  plane $e$. Now $e\cap e_i$ cuts $F$ in $n$ points, 
$i=1,3$; therefore the line in the configurations satisfying $p_1p_3g$ 
is determined by one of the $n^2$ pairs of points, $P_1$ on 
$e\cap e_1\cap F$ and $P_3$ on $e\cap e_3\cap F$; for such a choice of 
$P_1$ and $P_3$, we still can choose $P_2$ et $P_4$ among the 
$n-2$ remaining points on the line. Therefore there are
$$
n^2(n-2)(n-3)
$$
possible configurations.
In order to determine $G$, notice that for a given line $g$, 
there is a total of\newline 
$n(n-1)(n-2)(n-3)$ pairs of distinct points in $g\cap F$.
Finally, $p_1^3=0$. So we get~:
$$
2\e_{22}g_e=4n^2(n-2)(n-3)-3n(n-1)(n-2)(n-3)
$$
and therefore
$$
\e_{22}g_e=\frac12n(n-2)(n-3)(n+3)
\quad .
$$
\subsection{Third example: the number of lines tangent to a surface at four 
distinct points}
Denote by $\e_{2222}$  the condition that a line is tangent to a 
surface $F$ of degree $n$ at four distinct points; that amounts to say 
that among the $n$ points of intersection of the line with the 
surface, $4$ pairs come to coincide, say $(p_{1}, p_{2})$, $(p_{3},p_{4})$, 
$(p_{5},p_{6})$ and $(p_{7}, p_{8})$. It follows from the fact that 
the grassmannian ${\cal G}$ of lines in $\p$ is of dimension $4$ 
that for a generic surface there is a finite number of such 
quadritangent lines. We shall compute their number in a similar way 
as for the second example; it corresponds to formula $21)$ in 
\cite[pages 232 through 237]{Schubert:1879}. We will label some 
equations with letters, which have no analogue in \cite{Schubert:1879}.

It follows from the coincidence formula $(\spadesuit)$ that
$$
4!\e_{2222}=
(p_{1}+p_{2}-g)(p_{3}+p_{4}-g)
(p_{5}+p_{6}-g)(p_{7}+p_{8}-g)
$$
where the coefficient $4!$ is due to the fact that the roles of the 
$4$ pairs of points that come to coincide, $(p_{1}, p_{2})$, $(p_{3},p_{4})$, 
$(p_{5},p_{6})$ and $(p_{7}, p_{8})$, can be permuted.

The symbols $p_{i}p_{j}$, $i\not=j$ all express the same condition, 
therefore we deduce~:
\begin{equation*}
	\begin{split}
4!\e_{2,2,2,2}=&2^4p_{1}p_{2}p_{3}p_{4}-2^3\cdot 4g p_{1}p_{2}p_{3}
+2^2\binom{4}{2}g^2p_{1}p_{2}-2\binom{4}{3}g^3p_{1}+g^4\\
=&16p_{1}p_{2}p_{3}p_{4}-32g
p_{1}p_{2}p_{3}+24g^2p_{1}p_{2}-8g^3p_{1}+g^4
\end{split}\tag{a}
\end{equation*}

(Ther is a misprint in \cite[page 234, line -7]{Schubert:1879}~: 
$g_{4}$ is written instead of
$g^4$.)

It follows from formula I) that~:
\begin{equation*}
gp=p_{g}+g_{e}\implies
gp_{1}p_{2}p_{3}=p_{1}^2p_{2}p_{3}+g_{e}p_{1}p_{2}\tag{b}
\end{equation*}
Also~:
\begin{equation*}
g^2\,\stackrel{9)}{=}\, g_{p}+g_{e}\implies
g^2p_{1}p_{2}=\underbrace{g_{p}p_{1}p_{2}}_{=G}+g_{e}p_{1}p_{2}\tag{c}
\end{equation*}
and
\begin{equation*}
g^3p_{1}\,\stackrel{9)}{=}\, 
(g_{e}+g_{p})gp_{1}\,\stackrel{10)\hbox{ et }11)}{=}\,
2g_{s}p_{1}=2G\tag{d}
\end{equation*}
therefore we can use (b), (c) and (d) in (a) and find~:
\begin{equation*}
4!e_{2222}=16p_{1}p_{2}p_{3}p_{4}-32p_{1}^2p_{2}p_{3}-8g_{e}p_{1}p_{2}+10G\tag{e}
\end{equation*}
We compute separately each of the terms appearing in (e); we shall work in
$$
\big(F\tilde\times F\big)\times\big(F\tilde\times F\big)\times\big(F\tilde\times F\big)\times\big(F\tilde\times F\big)\times\G
$$
\begin{enumerate}
	\item\fbox{$G$}
	
	This represents the number of $8$-uples of distinct points
	$P_{1},\dots ,P_{8}$ that can be choosen in the intersection of a 
	generic line with the surface $F$ of degree $n$, that is~:
	$$
	G=n(n-1)(n-2)(n-3)(n-4)(n-5)(n-6)(n-7)
	$$
	\item
	\fbox{$g_{e}p_{1}p_{2}$}
	
	In the second example, working with $2$ pairs of points, we found 
	the formula $g_{e}p_{1}p_{3}=n^2(n-2)(n-3)$. 
	
	In the present case, we still have to choose $4$ points among the 
	remaining $n-4$ points in the intersection of a line with the 
	surface. Therefore we find~:
	$$
	g_{e}p_{1}p_{2}=n^2(n-2)(n-3)(n-4)(n-5)(n-6)(n-7)
	$$
	\item\fbox{$p_{1}^2p_{2}p_{3}$}
	
	Let's call $e_{1},e'_{1}$ the planes expressing condition $p_{1}^2$,
	so that $P_{1}$ will be among the $n$ points of the intersection 
	$F\cap e_{1}\cap e'_{1}$. 
	Let $e_{2}$ and $e_{3}$ be the planes expressing conditions 
	$p_{2}$ and $p_{3}$ respectively, and let $\ell$ be the line 
	containing $P_{1},\dots ,P_{8}$.
	
	The point $P_{2}$ must lye on the cone over the curve
	$F\cap e_{2}$ with vertex $P_{1}$ , that we shall denote by $C_{P_{1}}(F\cap
	e_{2})$; it is of degree $n$. 
	The point
	$P_{3}$ must lye on $C_{P_{1}}(F\cap e_{2})$ and on the curve
	$F\cap e_{3}$, which is of degree $n$. For a given $P_{1}$,
	there are $n^2$ possible choices for $\ell$, $P_{2}$ et $P_{3}$; 
	we must discard the choices corresponding to the 
	$n$ points of the intesection $e_{2}\cap
	e_{3}\cap F$, because otherwise we would have $P_{2}=P_{3}$. We 
	are left then with $n^2-n=n(n-1)$ solutions.
	
	Since there are $n$ possible choices for $P_{1}$,
	we find $n^2(n-1)$ possibilities for $P_{1},P_{2},P_{3}$. 
	We still have to choose
	$P_{4},\dots P_{8}$ among the remaining $(n-3)$ points on $\ell\cap
	F$; therefore we get~:
	$$
	p_{1}^2p_{2}p_{3}=n^2(n-1)(n-3)(n-4)(n-5)(n-6)(n-7)
	$$

	\item\fbox{$p_{1}p_{2}p_{3}p_{4}$}
	
	Suppose first that we only have the four points
	$P_{1},\dots ,P_{4}$. Let $e_{1},e_{2},e_{3},e_{4}$ be the planes
	expressing respectively conditions $p_{1},p_{2},p_{3},p_{4}$.
	Let's find the degree of the ruled surface $F'$ consisting of the 
	lines $\ell$ touching $e_{2}\cap F$, $e_{3}\cap F$,
	$e_{3}\cap F$; we have to compute the intersection of $F'$ with a 
	generic line, which amounts to compute $gp_{2}p_{3}p_{4}=(g_{e}+p_{2}^2)p_{3}p_{4}$.
	We know by the previous formulas that $g_{e}p_{2}p_{3}=n^2(n-2)$
	and $p_{2}^2p_{3}p_{4}=n^2(n-1)$, therefore 
	$F'$ is of degree $n^2(2n-3)$. 
	By B\'ezout's theorem\footnote{Schubert deduces B\'ezout's theorem from the 
	coincidence formula (\cite[page 46]{Schubert:1879}).},
	the intersection $F\cap F'\cap e_{1}$ consists of 
	$n^3(2n-3)$ points, among which, according to 3)
	above~:
	\begin{itemize}
		\item[$\bullet$] $n^2(n-1)$ are in $e_{2}$
		\item[$\bullet$] $n^2(n-1)$ are in $e_{3}$
		\item[$\bullet$] $n^2(n-1)$ are in $e_{4}$
	\end{itemize}
	and the remaining $n^3(2n-3)-3n^2(n-1)=n^2(2n^2-6n+3)$ constitute
	$p_{1}p_{2}p_{3}p_{4}$. If we take now in account the possible 
	choices for $P_{5},\dots ,P_{8}$ among the remaining $n-4$ 
	points, we get~:
	$$ 
	p_{1}p_{2}p_{3}p_{4}=n^2(2n^2-6n+3)(n-4)(n-5)(n-6)(n-7)
	$$
\end{enumerate}
\begin{remarque}
	Following \cite[page 236]{Schubert:1879}, let us note that for 
	$n=3$, if we count only the possible choices for $P_{1}$ through 
	$P_{4}$, we found the well known fact that there are $27$ lines 
	that cut a cuiuc surface in four four points, which implies that 
	the lines lye entirely in the cubic surface. This is a degenerate 
	case of quadrisecant lines; however, it can be shown that 
	there are no multiplicities to be taken in account, by using the 
	fact that along a line $\ell\subset F$, the tangent plane to $F$ 
	at a point $x\in \ell$ varies with $x$.
\end{remarque}
It remains to substitute the above formulas in (e)~:
\begin{multline*}
4!\e_{2222}=16n^2(2n^2-6n+3)(n-4)(n-5)(n-6)(n-7)-32n^2(n-1)(n-3)(n-4)(n-5)(n-6)(n-7)\\
-8n^2(n-2)(n-3)(n-4)(n-5)(n-6)(n-7)
+10=n(n-1)(n-2)(n-3)(n-4)(n-5)(n-6)(n-7)
\end{multline*}
and it follows that
$$
\e_{2222}=\frac{1}{12}n(n-4)(n-5)(n-6)(n-7)(n^3+6n^2+7n-30)
$$
which is the number of lines tangent to $F$ at four distinct points.

\bibliographystyle{plain}
\bibliography{eschubert}
\end{document}